\documentclass[letterpaper, 10pt, conference]{ieeeconf}      %

\IEEEoverridecommandlockouts                              %

\usepackage{graphicx} %
\usepackage{epsfig} %
\usepackage{amsmath} %
\usepackage{amssymb}  %
\usepackage{physics}
\usepackage{epstopdf}

\title{\LARGE \bf
Dynamic variable step size LMS adaptation algorithms---Application to adaptive feedforward noise attenuation*
}

\author{Tudor-Bogdan Airimitoaie$^{a}$, Bernard Vau$^{b}$, Dariusz Bismor$^{c}$, Gabriel Buche$^{d}$, and Ioan Doré Landau$^{d}$%
\thanks{*This work was not supported by any organization}%
\thanks{$^{a}$Tudor-Bogdan Airimitoaie is with the Univ. Bordeaux, CNRS, Bordeaux INP, IMS, 33405 Talence, France
        {\tt\small tudor-bogdan.airimitoaie@u-bordeaux.fr}}%
\thanks{$^{b}$Bernard Vau is with Exail, 12 avenue des Coquelicots, 94385 Bonneuil sur Marne, France
        {\tt\small bernard.vau@exail.com}}%
\thanks{$^{c}$Dariusz Bismor is with the Silesian University of Technology, ul. Akademicka 2A, 44-100 Gliwice, Poland
        {\tt\small dariusz.bismor@polsl.pl}}%
\thanks{$^{d}$Ioan Doré Landau and Gabriel Buche are with the Univ. Grenoble Alpes, CNRS, Grenoble INP, GIPSA-lab, 38000 Grenoble, France
        {\tt\small firstname.lastname@gipsa-lab.grenoble-inp.fr}}%
}

\newtheorem{thm}{Theorem}[section]

\newtheorem{lem}[thm]{Lemma}

\begin{document}

\maketitle
\thispagestyle{empty}
\pagestyle{empty}

\begin{abstract}
The paper explores in detail the use of dynamic adaptation gain/step size (DAG) for improving the adaptation  transient performance of variable step-size LMS (VS-LMS) adaptation algorithms. A generic form for the implementation of the DAG  within the VS-LMS algorithms is provided. Criteria for  the selection of the coefficients of the DAG filter which is required to be  a strictly positive real transfer operator are given.  
The potential of the VS-LMS adaptation algorithms using a DAG is then illustrated by experimental results obtained on a relevant adaptive active noise attenuation system. 
\end{abstract}

\section{INTRODUCTION}
The modern development of adaptation techniques in automatic control and signal processing started at the end of the fifties and beginning of the sixties (20th century). The paper \cite{Widrow60} introduced a gradient-based adaption algorithm in the discrete-time, later named the least mean squares (LMS). While the choice of the adaptation gain/step size for assuring the stability of the system was an open problem, interesting applications in the field of signal processing have been done. The paper \cite{Widrow75} gives an account of the applications of the LMS algorithm up to 1975.

In automatic control, the first attempt to synthesize adaptation algorithms has been probably the paper \cite{Whitaker58}, where a continuous-time formulation of a gradient type algorithm has been proposed. Unfortunately, dealing with feedback control systems to which a non-linear/time-varying loop (the adaptation loop) is added, raised crucial stability issues. The problem of the choice of the adaptation gain (step size) assuring the stability of the full system is fundamental.
Discrete time adaptation  algorithms assuring global asymptotic stability for any values of the adaptation gain were available since 1971 (\cite{LandauAsilomar71,LandauIFAC73_1,LandauIFAC75}). The concepts of "a priori" and "a posteriori" adaptation error emerged as key points for understanding the stability issues in the discrete-time context. 
These algorithms, derived from stability considerations, can be interpreted in the scalar case as gradient type algorithms trying to minimize a quadratic criterion in terms of the "a posteriori" prediction error (\cite{Landau86b,Landau90c}).

The signal processing community has concentrated its efforts in developing variable step size/adaptation gain algorithms in a scalar context (more exactly using a diagonal matrix adaptation gain) in order to improve the performance of the LMS algorithm. An exhaustive account of the various variable step-size LMS (VS-LMS) algorithms is provided by \cite{Bismor16}, where a unified presentation is done as well as an extensive comparison of the algorithms.
  It turns out that the adaptation algorithms developed in control from a stability point of view can also be interpreted as "variable step size" LMS algorithms.

In using adaptive/learning recursive algorithms there is an important problem to be addressed: the compromise between alertness (with respect to environment changes---like plant or disturbance characteristics) and stationary performances when using a constant value for the adaptation gain/step size. Accelerating the adaptation transient without augmenting the value of the adaptation gain/step size is a challenging problem.

Recently, the concept of dynamic adaptation gain (DAG) has been introduced in  \cite{LandauCDC22,LandauAuto23,LandauJSV23} as means to accelerate significantly the adaptation transients without modifying the steady state (asymptotic) properties of an algorithm for a given adaptation gain/step size. The correcting term in the adaptation algorithm is first filtered before its use for the estimation of a new parameter value. With an appropriate choice of the parameters of this filter, which should be characterized by a strictly positive real (SPR) transfer function, a significant improvement of the adaptation transient is obtained. The design of this filter is well understood and design tools are available.

The main objective of this paper is to show that the DAG introduced in the context of stability based adaptation algorithms can be successfully applied to VS-LMS adaptation algorithms leading to similar significant acceleration of the adaptation transients. This will be illustrated by real-time experiments on an adaptive feedforward attenuator (a silencer).
\section{Introducing the DAG-VS-LMS algorithms}\label{DAG}
\begin{figure}[!htb]
    \begin{center}
    \includegraphics[width=0.6\columnwidth]{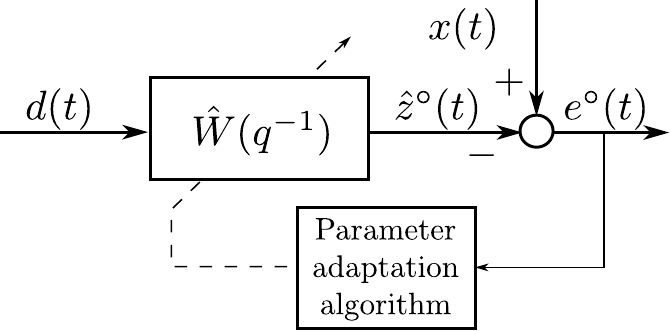}
    \caption{Least mean squares (LMS) adaptive filtering problem.}
    \label{fig_LMS_basic_scheme}
    \end{center}
\end{figure}
The VS-LMS (variable step size least mean squares) algorithms are very popular in the field of signal processing and the field of active vibration and noise control. There is a strong similarity with some of the algorithms used in the adaptive control and recursive system identification. The VS-LMS algorithms (which are improvements of the original LMS algorithm) will be briefly reviewed in order to add the dynamic adaptation gain/step size introduced in \cite{LandauCDC22}. The aim of the LMS parameter adaptation/learning algorithm (PALA) is to drive the parameters of an adjustable model in order to  minimize a quadratic criterion in terms of the prediction error (difference between real data and the output of the model used for prediction).

The basic block diagram illustrating the LMS algorithm's operation is shown in Fig.~\ref{fig_LMS_basic_scheme}. The adaptive filter $W(q^{-1})$ is fed with the input sequence
 $d(t)$. The output of the filter, $z(t)$, is compared with the desired signal, $x(t)$, to compute the error signal $e(t)$. The LMS algorithm adjusts the weights of the $W(q^{-1})$ filter to minimize the error. 

Consider that the desired signal can be described by:
\begin{align}
  x(t)=\theta^T\phi(t),
  \label{eq:output}
\end{align}
where the \textit{parameter vector} and the \textit{measurement vector} are denoted by
\begin{align}
  \theta^T&=\left[ w_0,w_1,\ldots,w_{n_W} \right],\text{ and}\\
  \phi^T(t)&=\left[ d(t),d(t-1),\ldots,d(t-n_W) \right],
\end{align}
respectively. The adjustable prediction model of the adaptive filter will be described by:
\begin{align}
  \hat{z}^\circ(t)=\hat{\theta}^T(t-1)\phi(t),
\end{align}
where $\hat{z}^\circ(t)$ is termed the \textit{a priori} predicted output depending upon the values of the estimated parameter vector $\theta$ at instant $t-1$:
\begin{align}
  \hat{\theta}^T(t-1)&=\left[ \hat{w}_0(t-1),\hat{w}_1(t-1),\ldots,\hat{w}_{n_W}(t-1) \right].
\end{align}
It is very useful to consider also the \emph{a posteriori} predicted output computed on the basis of the new estimated parameter vector at $t$, $\hat{\theta} (t)$, which will be available somewhere between $t$ and $t+1$. The \emph{a posteriori} predicted output will be given by:
\begin{align}
\hat{z}(t)&=\hat{\theta}^T(t) \phi(t)
\label{a postout}
\end{align}
One defines an \emph{a priori} prediction error as:
\begin{equation}
 e^\circ(t)=x(t)-\hat{z}^\circ(t)= [\theta - \hat{\theta}(t-1)]^T\phi(t) \label{eq_apriori}
\end{equation}
and an \emph{a posteriori} prediction error\index{a posteriori prediction error} as:
\begin{equation}
e(t)=x(t)-\hat{z}(t)= [\theta - \hat{\theta}(t)]^T\phi(t). \label{eq_aposteriori}
\end{equation}
The VS-LMS algorithms update the filter taps according to the formula (see \cite{Widrow75}):
\begin{align}
  \hat{\theta}(t)=\hat{\theta}(t-1)+\mu(t) \phi(t) e^\circ(t),
  \label{eq_vslms}
\end{align}
where $\mu(t)$ is the variable step-size parameter. In the standard form of the LMS algorithm, the step-size has a constant value $\mu(t)=\mu$. Large values of $\mu$ allow for fast adaptation, but also give large excess mean square error (EMSE, see \cite{Bismor16}). Too large step-sizes may lead to the loss of stability (see \cite {LandauSpringer11} for a discussion). On the other hand, too small step-sizes give slow convergence, which in many practical applications is not desirable. A variable step-size $\mu(t)$ can provide a compromise. The first VS-LMS algorithm was the Normalized LMS (NLMS) algorithm proposed in 1967 independently by \cite{Albert67,Nagumo67}, which uses the following equation for the step-size:
\begin{align}
  \mu(t)=\frac{\mu}{\delta+\phi(t)^T\phi(t)},
\end{align}
where $\delta$ is a very small scalar value used in order to avoid division by zero (typically of the order of $10^{-16}$).
\footnote{The case of NLMS using larger values of $\delta$ (for example $\delta=1$) is discussed in \cite{Goodwin84a}.}

In adaptive control, where the stability of the full adaptive control system is considered as a fundamental issue, VS-LMS algorithms have been developed from a stability point of view. However, when scalar type adaptation gain is used (i.e. diagonal matrix adaptation gain), these algorithms  can be interpreted as gradient type algorithms for the minimization of a quadratic error in terms of the \textit{a posteriori} prediction error defined in \eqref{eq_aposteriori}. See \cite{LandauSpringer11} for details. See also \cite{LandauJSV23} for a derivation of this type of algorithm and its application to active noise control. This algorithm will be termed PLMS  (to distinguish it with respect to the LMS and the other VS-LMS which are using the \textit{a priori }prediction error). The algorithm has the form:
\begin{align}
  \hat{\theta}(t)&=\hat{\theta}(t-1)+\mu \phi(t) e(t)\\
  &=\hat{\theta}(t-1)+\mu(t) \phi(t) e^\circ(t),
  \label{eq_plms}
\end{align}
where:
 \begin{align}
  \mu(t)=\frac{\mu}{1+\mu\phi(t)^T\phi(t)},
  \label{eq_mu}
\end{align}
See Appendix~\ref{PLMS_pr} for a direct derivation of the algorithm.

When using the \textit{dynamic adaptation gain/learning rate} (DAG)  equation \eqref{eq_vslms} of the VS-LMS algorithms will take the form:
\begin{equation}
\hat{\theta}(t)=\hat{\theta}(t-1)+\frac{C(q^{-1})}{D'(q^{-1})}\left[ \mu(t) \phi(t) e^\circ(n) \right]
\label{eq:dynamic1}
\end{equation}
where\footnote{The complex variable $z^{-1}$ will be used for characterizing the system's behaviour in the frequency domain and the delay operator $q^{-1}$ will be used for describing the system's behavior in the time domain.} $\frac{C(q^{-1})}{D'(q^{-1})}$ is termed the ``dynamic adaptation gain/learning rate" (DAG) and has the form:
\begin{equation}\label{eq_theta_DAG}
H_{DAG}(q^{-1})=\frac{C(q^{-1})}{D'(q^{-1})}=\frac{1+c_{1}q^{-1} + \ldots+c_{n_C}q^{-n_C}}{1-d'_{1}q^{-1}-\ldots -d'_{n_{D'}}q^{-n_{D'}}}
\end{equation}
The effective implementation of the algorithm given in \eqref{eq:dynamic1} leads to:
\begin{multline}
 \hat{\theta} (t) =d_1\hat{\theta}(t-1)+ d_2\hat{\theta}(t-2)+\ldots+d_{n_D}\hat{\theta}(t-n_D) \\
 +\mu(t)\phi(t)e^\circ(t)+ c_1\mu(t-1)\phi(t-1)e^\circ(t-1)+ \\
 + \ldots +c_{n_C}\mu(t-n_C)\phi(t-n_C)e^\circ(t-n_C)
 \label{eq:ARMA3}
\end{multline}
where ($n_D=n_{D'}+1$):
\begin{equation}
d_i=(d'_i-d'_{i-1})~~; i=1,... n_D; d'_0=-1,~d'_{n_D}=0
\label{ARIMA2}
\end{equation}
To implement the algorithm, one needs a computable expression for $e^\circ(n)$. This is obtained by computing $\hat{z}^\circ(n)$ in \eqref{eq_apriori} as:
\begin{align}
  \hat{z}^\circ(t)=\hat{\theta}_0^T(t-1)\phi(t),
  \label{eq_z0DAG}
\end{align}
where
\begin{multline}
  \hat{\theta}_0^T(t-1)=d_1\hat{\theta}(t-1)+ d_2\hat{\theta}(t-2)+\ldots+d_{n_D}\hat{\theta}(t-n_D) \\
  +c_1\mu(t-1)\phi(t-1)e^\circ(t-1)+ \\
  + \ldots +c_{n_C}\mu(t-n_C)\phi(t-n_C)e^\circ(t-n_C).
\end{multline}

\subsection{Relations with other algorithms}
Many algorithms have been proposed for accelerating the speed of convergence of the adaptation algorithms derived using the "gradient rule". The algorithm of \eqref{eq:dynamic1} is termed ARIMA (Autoregressive with Integrator Moving Average). As discussed in \cite[Section 8]{LandauAuto23}, a number of well known algorithms are particular cases of the ARIMA algorithm. The various algorithms described in the literature are of MAI form or ARI form. The MAI form includes ``Integral+ Proportional" algorithm \cite{LandauSpringer11,AirimitoaieAuto13} ($c_1\neq  0, c_i=0,\forall~i>1$, $d'_i=0,\forall~i>0$),
``Averaged gradient" ($c_i, i=1,2,...$, $d'_i=0,\forall~i>0$) \cite{Schmidt18,Pouyanfar17}. The ARI form includes ``Conjugate gradient" and ``Nesterov" algorithms \cite{Livieris14,AnnaswamyCDC19} ($c_i=0, i=1,2,..,d'_1 \neq  0, d'_i=0,i>1)$ as well as the ``Momentum back propagation" algorithm \cite{Jacobs88} which corresponds to the conjugate gradient plus a normalization of $\mu$ by $(1-d'_1)$\footnote{There are very few indications how to choose the various weights in the above mentioned algorithms.}.
A particular form of the ARIMA algorithms termed ``ARIMA2" ($c_1,c_2, c_i=0,\forall~i>2,d'_1\neq 0,d'_i=0,\forall~i>1$) will be studied subsequently and evaluated experimentally.\footnote{The algorithms mentioned above can be viewed as particular cases of the ARIMA2 algorithm.}
\section{Design of the dynamic adaptation gain/learning rate}\label{design}
\subsection{Performance issues}
The dynamic adaptation gain/learning rate will introduce a phase distortion on the gradient depending upon the frequency. Assume that the algorithms should operate for all frequencies in the range: $0$ to $0.5f_s$ ($f_s$ is the sampling frequency). 
Assume that the gradient of the criterion to be minimized contains a single frequency.
 In order to minimize the criterion, the phase distortion introduced by the dynamic adaptation gain/learning rate should be less than $90^\circ $ at all the frequencies. 
In other terms, the transfer function $\frac{C(z^{-1})}{D'(z^{-1})}$ should be strictly positive real (SPR) (see \cite{LandauSpringer11} for a detailed definition).
 In order that a transfer function be strictly positive real, it must first have its zeros and poles inside the unit circle. One has the following property:
\begin{lem}\label{lem_I_equal_0}
  Assume that the polynomials $C(z^{-1})$ and $D'(z^{-1})$ have all their zeros inside the unit circle, then:
  \begin{equation}
  I=\int_{0}^{\pi}\log\left( \left|\frac{C(e^{-i\omega})}{D^{'}(e^{-i\omega} )}\right | \right)\dd  \omega =0.
  \label{eq:integral}
  \end{equation}
\end{lem}
The proof relies on the Cauchy Integral formula (see \cite{LandauAuto23}).

\begin{figure}[htb]
  \includegraphics[width=1\columnwidth]{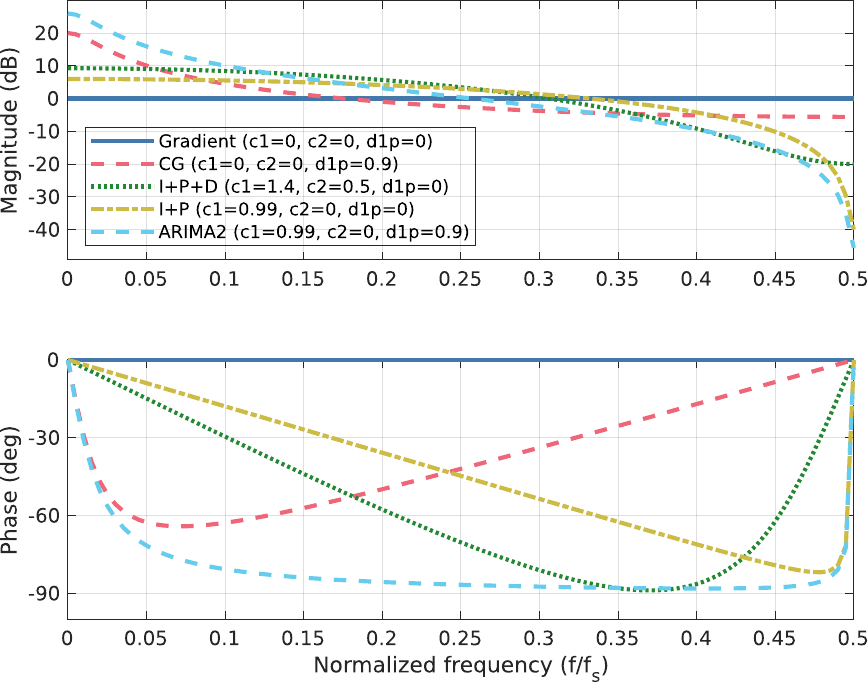}
  \caption{Frequency characteristics of DAGs used in the experimental section (see also Table~\ref{perf1}).}
  \label{fig_bode_comp}
\end{figure}

This result allows to conclude that the average gain over the frequency range $0$ to $0.5f_s$ is 0~dB, i.e. on the average this filter will not modify the adaptation gain/step size. It is just a frequency weighting of the adaptation gain/step size. To be more specific, Figure~\ref{fig_bode_comp}  shows the frequency characteristics of various DAG's that will be subsequently used in the experimental section (the VS-LMS algorithm corresponds to the axis $0$~dB).
 It can be observed first that the phase is within the range $\pm 90^\circ$, i.e. they are SPR. Then one can observe that effectively the average gain over the frequency range $0$ to  $0.5f_s$ ($f_s=2500$ Hz) is $0$~dB. Now examining the magnitude, one observes that all are low pass filters amplifying low frequencies. This means that if the frequency content of the gradient is in the low frequency range, the effective adaptation gain/learning rate will be larger than $\mu$, which should have a positive effect upon the adaptation/learning transient. In particular, the DAG which has a larger gain in low frequencies (ARIMA2) should provide the best performance  (which is indeed the case---see Section \ref{Exp}).

Since we need to have a DAG which is SPR,  we will provide subsequently the tools for design of a SPR DAG. We will consider the case of the ARIMA2 algorithm introduced in \cite{LandauCDC22}. The DAG in this case will have the form:
\begin{align}\label{eq_H_DAG}
H_{DAG}=\frac{C(q^{-1})}{D'(q^{-1})}=\frac{1+c_1q^{-1}+c_2q^{-2}}{1-d'_1 q^{-1}}
\end{align}
A criterion for the selection of $c_1$, $c_2$ and $d'_1$ in order that the DAG be SPR is given below.
\begin{lem}\label{lemma_SPR}
The conditions assuring that $ H_{DAG}(z)=\frac{1+c_1z^{-1}+c_2z^{-2}}{1-d_1^{'}z^{-1}}$ is strictly positive real (SPR) are:
\begin{itemize}
\item for $c_2\leq 0$, $c_1$ must be such that
\[-1 -c_2 < c_1 < 1+c_2 \]

\item for  $c_2\geq 0$
\begin{itemize}
\item if the following condition is satisfied \[ 2(d_1^{'}-c_2) < \sqrt{2(c_2-c_2^2)(1-d_1^{'2})} <  2(d_1^{'}+c_2)  \]
  the maximum bound on $c_1$ is given by  \[c_1 < d_1^{'}-3d_1^{'}c_2+2\sqrt{2(c_2-c_2^2)(1-d_1^{'2})} \]
	otherwise the maximum bound on $c_1$ is given by  $c_1 < 1+c_2$ 
\item   if the following condition is satisfied \[ 2(d_1^{'}-c_2) < -\sqrt{2(c_2-c_2^2)(1-d_1^{'2})} <  2(d_1^{'}+c_2)  \]
the minimum bound on $c_1$ is given by \[c_1 > d_1^{'}-3d_1^{'}c_2-\sqrt{2(c_2-c_2^2)(1-d_1^{'2})} \]
	otherwise the minimum bound on $c_1$ is given by  $c_1 > -1-c_2$ 
\end{itemize}

\end{itemize}

\end{lem}
The proof of this result is given in \cite{LandauAuto23}.

From the conditions of Lemma \ref{lemma_SPR}, closed contours in the plane $c_2-c_1$ can be defined for the different values of $d'_1$ allowing to select $c_1$ and $c_2$ for a given value of $d'_1$ such that the DAG be SPR. 
Note that a necessary condition for the selection of the parameters $c_1, c_2, d'_1$ is that both the denominator and the numerator of the $H_{DAG}$ should be asymptotically stable.
\subsection{Stability issues}
If one wants to use VS-LMS algorithms with large values of the adaptation gain $\mu$, the stability analysis of the resulting scheme using a dynamic adaptation gain is an important issue. For the case of the PLMS algorthm, this analysis has been carried in detail (see \cite{LandauAuto23, LandauJSV23}). We will just indicate the main result subsequently.
Equation \eqref{eq:dynamic1} can be expressed also as:
\begin{equation}
(1-q^{-1})\hat{\theta}(t+1)=+\frac{C(q^{-1})}{D'(q^{-1})}\left[ \mu(t) \phi(t) e^\circ(t+1) \right] \nonumber
\label{eq:dynamic1''}
\end{equation}
leading to:
\begin{equation}
\hat{\theta}(t+1)=H_{PAA}(q^{-1})[\mu\phi(t) e(t+1)],
\label{eq:genalg3}
\end{equation}
where  $e(t+1)=e^\circ(t+1)(1+\mu\phi^T(t)\phi(t))^{-1}$ and $H_{PAA}=(1-q^{-1})^{-1}H_{DAG}$ is a MIMO diagonal transfer operator having identical terms. All the diagonal terms are described by:
\begin{align}
H_{PAA}^{ii}(q^{-1})&=\frac{1+c_{1}q^{-1} +\ldots+c_{n_C}q^{-n_C}}{(1-q^{-1})(1-d'_{1}q^{-1}-\ldots-d'_{n_{D'}}q^{-n_D'})}\nonumber \\
&=\frac{C(q^{-1})}{(1-q^{-1})D'(q^{-1})}=\frac{C(q^{-1})}{D(q^{-1})}.
\label{eq:ARIMA1}
\end{align}
The relation between the coefficients of polynomials $D$ and $D'$ is given  in \eqref{ARIMA2}.

One has the following result:
\begin{thm}\label{thm_PR}
  For the system described by Equations~\eqref{eq:output} through \eqref{eq_aposteriori} using the PLMS  algorithm of \eqref{eq:dynamic1} and \eqref{eq_theta_DAG}
  one has $\lim_{t \to \infty} e(t+1) =0$ for any positive adaptation gain $\mu$ and any initial conditions $\theta (0), e (0)$, if $H_{PAA}^{ii}(z^{-1})$  given in \eqref{eq:ARIMA1} is a PR transfer function with a pole at $z=1$.
\end{thm}

The proof is given in \cite{LandauAuto23}. %

The adaptive/learning system considered in the Theorem~\ref{thm_PR} leads to an equivalent feedback representation where the equivalent feedforward path is a constant positive gain and the equivalent feedback path features the $H_{PAA}$ (see \cite{LandauCDC22}). However, in many cases the equivalent feedforward path will be a transfer operator. In such situations in addition to the  PR condition upon the $H_{PAA}$, there will be an additional SPR condition upon the transfer operator characterizing the equivalent feedforward path.

For small values of the adaptation gains/learning rates, the passivity/stability condition can be relaxed using \textit{averaging}  \cite{Anderson86}.
  In this case, the most important is that the $H^{ii}_{PAA}$ is PR in the frequency region of operation (mainly defined by the spectrum of the input signals to the systems).

\begin{figure}[htb!]%
    \begin{center}
    \includegraphics[width=\columnwidth]{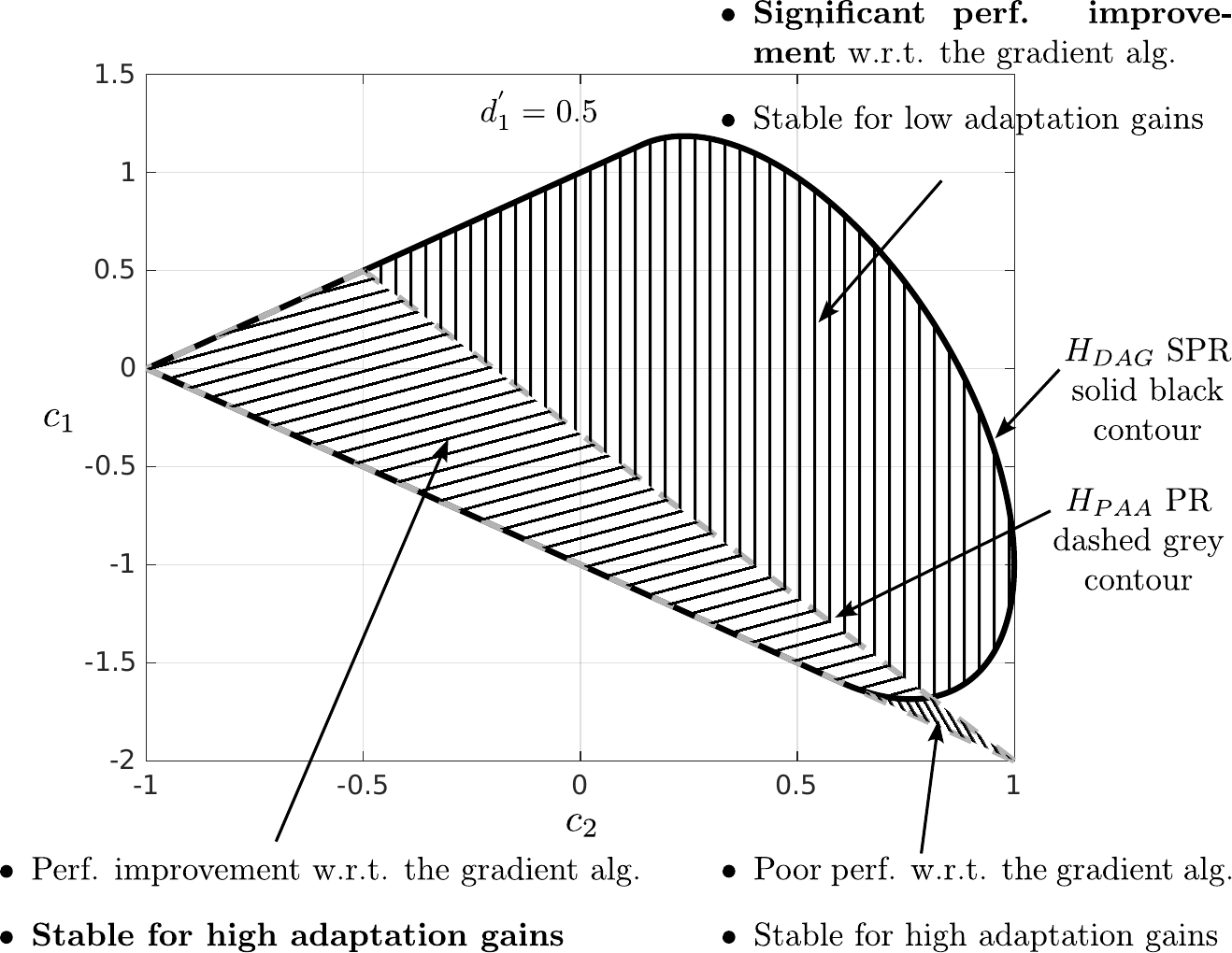}
    \caption{Intersection in the plane $c_1 - c_2$ of the contour  $H_{PAA}=PR$ with the contour  $H_{DAG}=SPR$ for $d'_1=0.5$.}
    \label{fig_GADww}
    \end{center}
\end{figure}

It is interesting to see intersections of the contours assuring the SPR of the $H^{ii}_{DAG}$ with the contours assuring that $H^{ii}_{PAA}$ is PR. 
Such an intersection is shown in Fig.\ref{fig_GADww} for the case of the $H_{DAG}$ given in \eqref{eq_H_DAG}. From this figure one can conclude that not all the SPR $H_{DAG}$ will lead to a PR $H_{PAA}$. In such cases, the performance is improved for low adaptation gains, but one can not guarantee asymptotic stability for large values of the adaptation gain.  Fig.~\ref{fig_GADww} shows also that there is a region where despite  that $H_{PAA}$ is PR, $H_{DAG}$ is not SPR. For such configurations, large adaptation gains can be used but the adaptation transient is slower than for the basic gradient algorithm.

\section{Experimental results}\label{Exp}
\begin{figure}[htb!]
  \centering
  
		\includegraphics[width=\columnwidth]{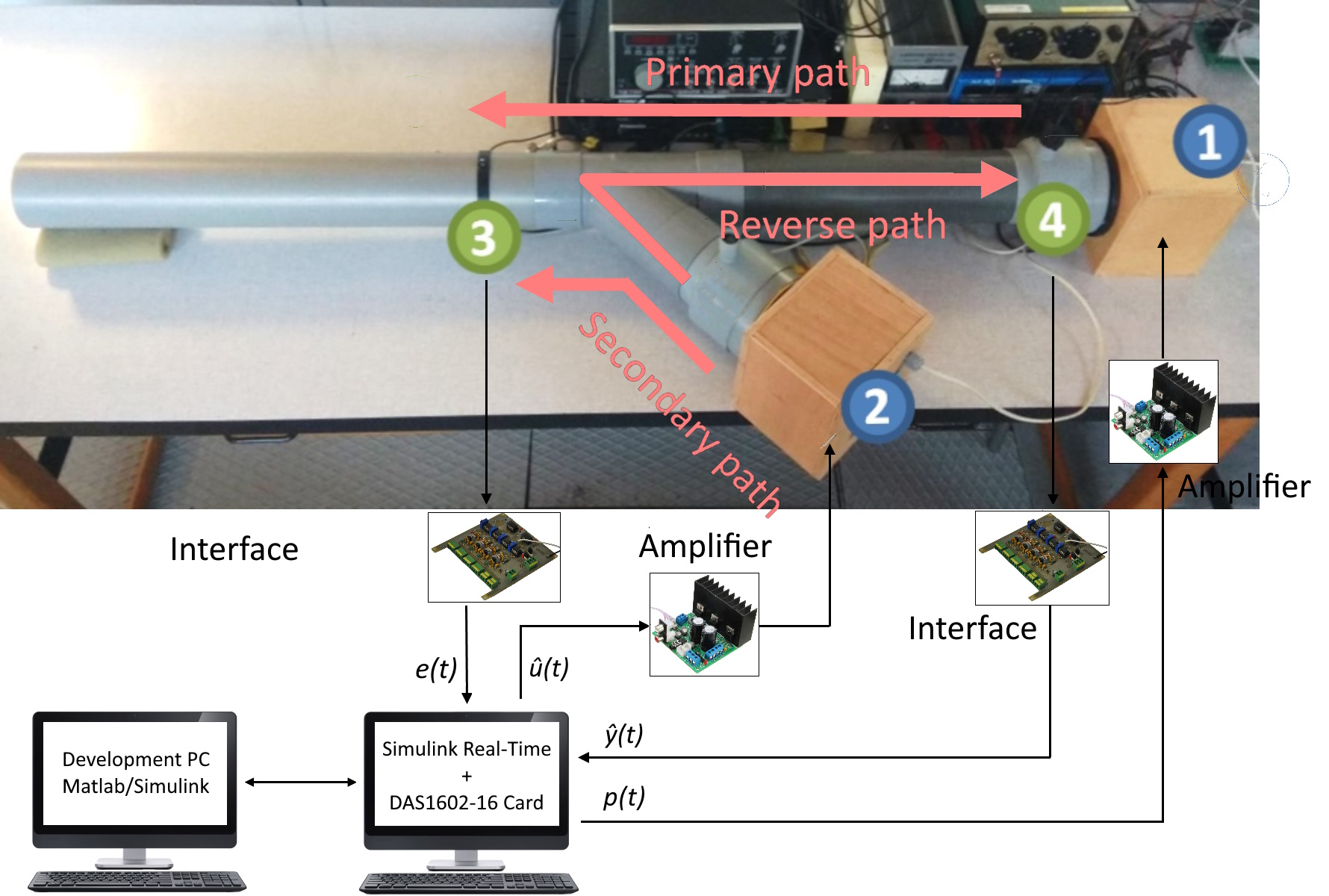} %
	\caption{Duct active noise control test-bench .}
	\label{photos_test_bench}
\end{figure}

The objective of his section is to show that the dynamic adaptation gain can be implemented on various VS-LMS algorithms and this will lead to a significant acceleration of the adaptation transient. Specifically for this paper the LMS, NLMS and PLMS algorithms have been implemented and tested experimentally on an active noise control test-bench (adaptive feedforward noise attenuation).
 Figure~\ref{photos_test_bench} shows the view of the test-bench used for experiments. Detailed description can be found in \cite{LandauAuto23}. 
The speakers and the microphones were connected to a target computer with the Simulink Real-Time\textsuperscript{\textregistered}. A second computer is used for development and operation with Matlab/Simulink. The sampling frequency was $f_s=2500$~Hz. 

The models of the various paths  in Figure \ref{photos_test_bench} are characterized by the presence of many pairs of very low damped poles and zeros. These models have been identified experimentally. The  orders of the various identified models are: $n_G=33$, $n_M=27$ and $n_D=27$.

\begin{figure}[ht]%
    \begin{center}
    \includegraphics[width=8cm]{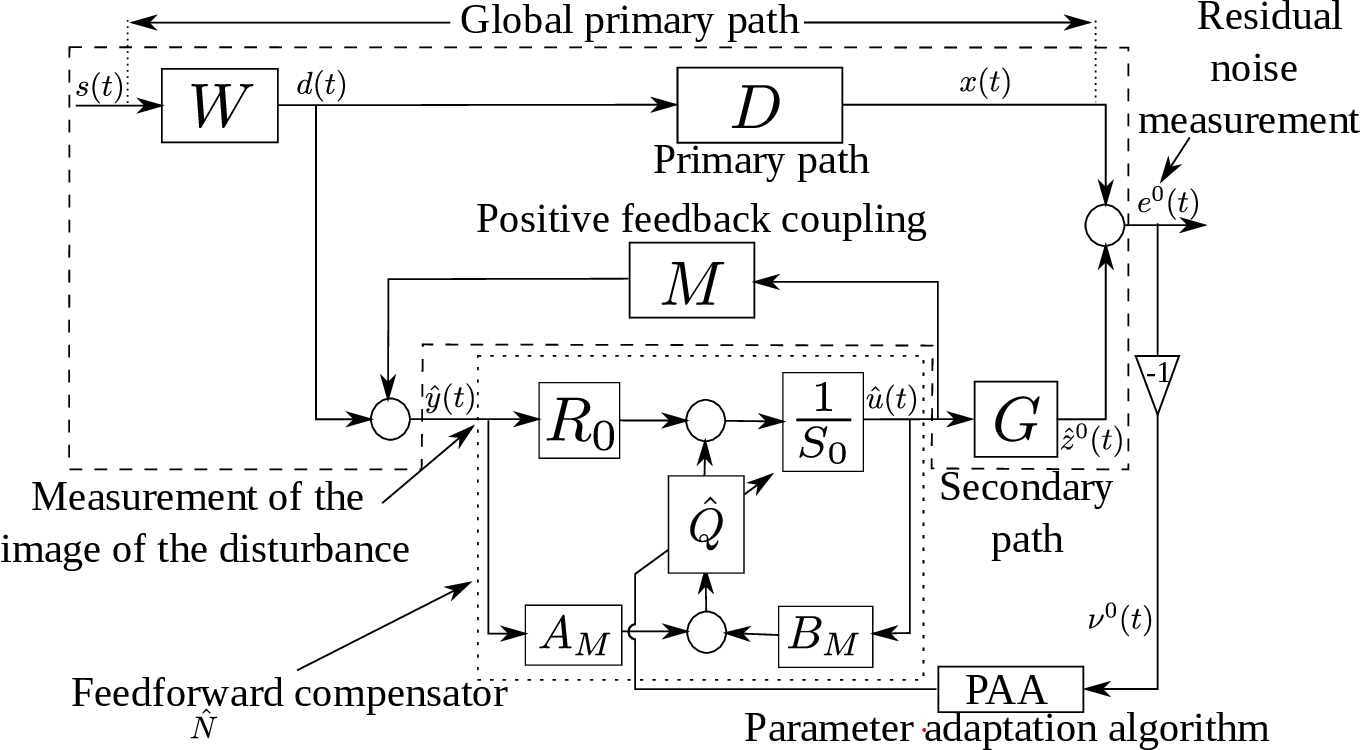}
    \caption{Feedforward AVC with FIR-YK adaptive feedforward compensator.}
    \label{fig_feedforwardAVC}
    \end{center}
\end{figure}

The objective is to attenuate an incoming unknown broad-band noise disturbance. The corresponding block diagram for the adaptive feedforward noise compensation using FIR Youla-Kucera (FIR-YK) parametrization of the feedforward compensator (introduced in \cite{LandauAuto12} for active vibration control and in \cite{airimitoaie:hal-02947816} for active noise control) is shown in Figure~\ref{fig_feedforwardAVC}.

The adjustable filter $\hat{Q}$ has the FIR structure:
\begin{equation}\label{eq_Qhat(q-1)}
            \hat{Q}(q^{-1}) = {\hat{q}_{0}}+\hat{q}_{1}q^{-1}+...+\hat{q}_{n_{Q}}q^{-n_{Q}}
        \end{equation}
and the parameters $\hat{q}_{i}$ will be adapted in order to minimize the residual noise.
The algorithm
  can be summarized as follows.
One defines:
\begin{eqnarray}
\theta^T&=&[q_0,q_1,q_2,\dots,q_{n_Q}]\label{teta}\\
		 \hat{\theta}^T&=&[\hat{q}_0,\hat{q}_1,\hat{q}_2,\dots,\hat{q}_{n_Q}]\label{tetahat}\\
		 \phi^T(t)&=&[v(t+1),v(t),\dots,v(t-{n_Q}+1)]\label{fi}
	\end{eqnarray}
where:
\begin{equation}
		 v(t+1)=B_{M}\hat{y}(t+1)-A_{M}\hat{u}(t+1)=B^{*}_{M}\hat{y}(t)-A_{M}\hat{u}(t+1)\label{eq:alpha}
	\end{equation}
	One defines also the regressor vector (a filtered observation vector) as:
 \begin{equation}\label{eq_phif}
\phi_{f}(t)=L(q^{-1})\phi(t)=[v_f(t+1),v_f(t),\dots,v_f(t-n_Q+1)]
\end{equation}
where
\begin{equation}
v_{f}(t+1)=L(q^{-1})v(t+1)\label{eq:filtalg2}
\end{equation}
Using $R_0=0$ and $S_0=1$, the poles of the internal positive closed loop are defined by $A_M$ and they will remain unchanged. The filter  used in \eqref{eq:filtalg2} becomes $L=\hat{G}$ and the associated linear transfer operator appearing in the equivalent feedforward path is
\begin{equation}\label{eq_Halgo2}
H(q^{-1})=\frac{G(q^{-1})}{\hat{G}(q^{-1})}
\end{equation}
(the algorithm uses an approximate gradient). The transfer function associated to $H(q^{-1})$ should be SPR in order to assure asymptotic stability in the case of perfect matching. This is a very mild condition as long as a good experimental identification of the models is done.

The VS-LMS algorithms which were used are of the form given in \eqref{eq_vslms} where $e^\circ(n)$ is given by \eqref{eq_apriori} with $z^\circ$ given in \eqref{eq_z0DAG}, where $\hat{\theta}$ is given by \eqref{tetahat} and $\phi$ is replaced by $\phi_f$ given in \eqref{eq_phif}\footnote{In signal processing literature when using a filtering of the regressor, the algorithms are termed FU-VSLMS.}.  The adjustable filter $\hat{Q}$ had 60 parameters.

The values of the adaptation gain $\mu$ for the three algorithms were tuned such that in the absence of the DAG, the performance was close for the three algorithms. A specificity of this application is also the low value of the average of $\phi(n)^T\phi(n) ~~(<<1)$. This means that the LMS and PLMS show a very close behavior for a given $\mu$ and that the adaptation gain $\mu$ for the NLMS should be much lower in order to get similar performances. The choices of $\mu$ are 0.2 for LMS, 0.0002 for NLMS and 0.22 for the PLMS.

\renewcommand{\tabcolsep}{4pt}
 \begin{table}[htb!]
    \begin{center}%
    \caption{Parameters of ARIMA2 Dynamic Adaptation Gain.}
        \begin{tabular}{c|c|c|c|c|c}\hline
            Algorithm &$H_{PAA}$--PR &$H_{DAG}$--SPR&$c_1$  		& $c_2$ & $d'_1$ 
            \\ \hline \hline
            Integral (gradient) & Y &Y&$0$ 	& $0$ 	& $0$ 	
            \\ \hline
            Conj.Gr/Nest..& N&Y& $0$ 	& $0$ 	& $0.9$ 	
            \\ \hline
            I+P+D & N &Y& $1.4$ &  $0.5$ & $0$  
            \\ \hline
            I+P & Y&Y & $0.99$ &  $0$ & $0$ 
            \\ \hline
            ARIMA 2 & N&Y& $0.99$ &  $0$ & $0.9$ 
            \\ \hline
        \end{tabular}
        \label{perf1}
    \end{center}  
\end{table}

\begin{figure}[htb]
    \includegraphics[width=1\columnwidth]{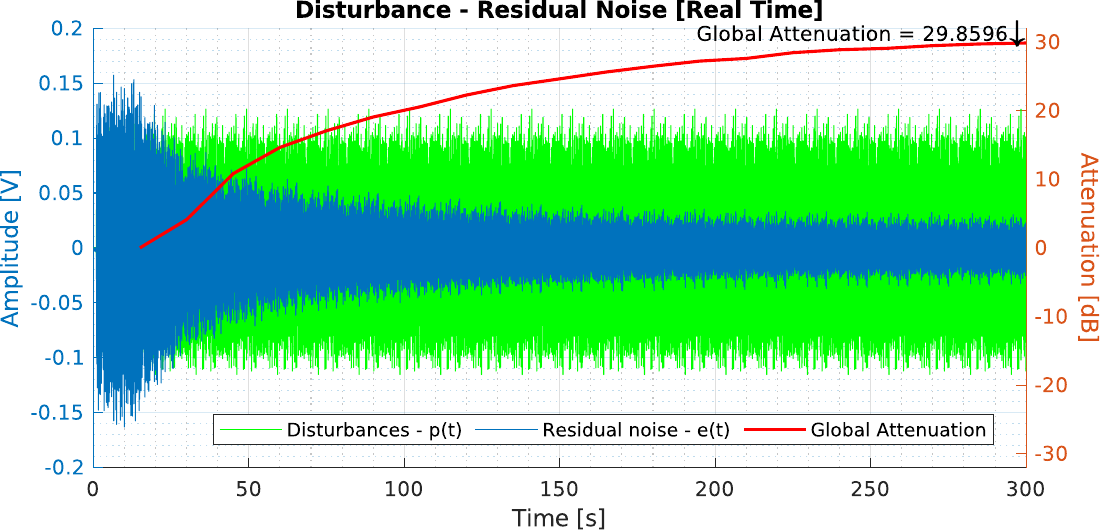}
		\includegraphics[width=1\columnwidth]{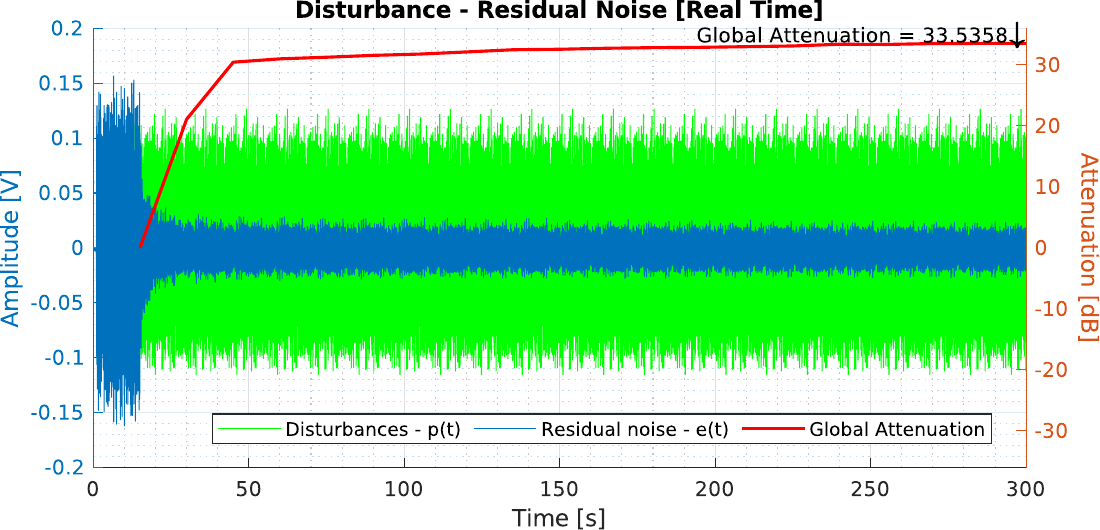}
	\caption{Time evolution of the residual noise using the standard NLMS adaptation algorithm (top) and using the NLMS + ARIMA2 algorithm (bottom), for $\mu=0.0002$.}
	\label{fig:resnoise}
\end{figure}

A broad-band disturbance 70--170 Hz was used as an unknown disturbance acting on the system. The steady state and transient attenuation\footnote{The attenuation is defined as the ratio between the variance of the residual noise in the absence of the control and the variance of the residual noise in the presence of the adaptive feedforward compensation. The variance is evaluated over an horizon of 3 seconds.} were evaluated for the various values of the parameters $c_1,~c_2$ and $d'_1$ given in Table~\ref{perf1}. The system was operated in open-loop during the first 15~s. Figure~\ref{fig:resnoise} shows the time response of the system as well as the evolution of the global attenuation when using the standard NLMS algorithm (top) and when the ARIMA2 DAG is incorporated (bottom) with  $c_1=0.99,~c_2=0,~d'_1=0.9$ (last row of Table~\ref{perf1}). One observes a significant acceleration of the adaptation transient. The acceleration obtained is equivalent to that obtained with an adaptation gain 25 times higher using the standard NLMS algorithms.
  Similar results were obtained for the LMS and the PLMS algorithms.

\begin{figure}[!htb]
  \begin{center}
  \includegraphics[width=\columnwidth]{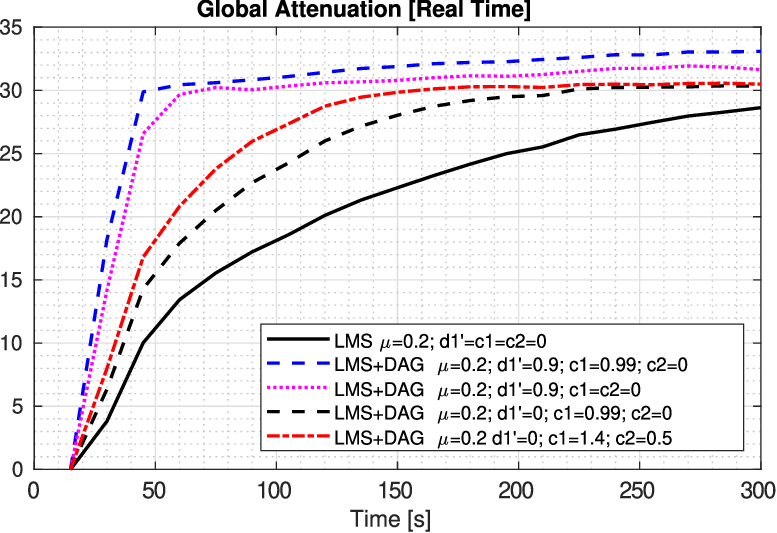}
  \caption{Time evolution of the global attenuation for the LMS algorithm with various DAGs, using $\mu=0.2$.}
  \label{fig_LMS_dag}
  \end{center}
\end{figure}

\begin{figure}[!htb]
  \begin{center}
  \includegraphics[width=\columnwidth]{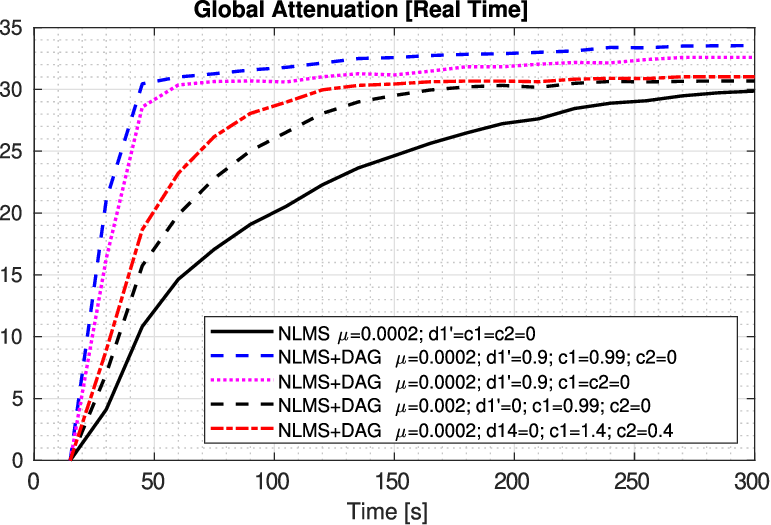}
  \caption{Time evolution of the global attenuation for the NLMS algorithm with various DAGs, using $\mu=0.0002$.}
  \label{fig_NLMS_dag}
  \end{center}
\end{figure}

\begin{figure}[!htb]
  \begin{center}
  \includegraphics[width=\columnwidth]{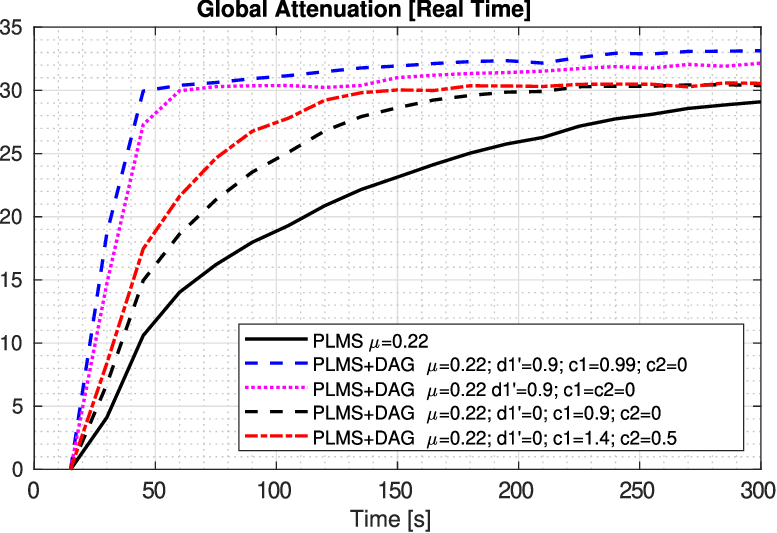}
  \caption{Time evolution of the global attenuation for the PLMS algorithm with various DAGs, using $\mu=0.22$.}
  \label{fig_PLMS_dag}
  \end{center}
\end{figure}

Figures \ref{fig_LMS_dag}, \ref{fig_NLMS_dag}, and \ref{fig_PLMS_dag} give the time evolution of the attenuation for the LMS, NLMS and PLMS algorithms and for the various DAG given in Table \ref{perf1}. As one can observe, the effect of the DAG is similar for the three algorithms.

\section{Conclusion}
\label{Conclusions} 
The paper emphasizes the potential of the dynamic adaptation/step size (DAG) for improving the adaptation transients  of VS-LMS  adaptation algorithms.  The main point is that the DAG should be characterized by a \textit{strictly positive real} (SPR) transfer function if we would like to operate correctly for any frequency in the range $0$ to $0.5$ of the sampling frequency.
  Experimental results on a relevant adaptive active noise control system have illustrated the feasibility and the performance improvement achieved using a DAG on VS-LMS adaptation algorithms.
\begin{appendices}
\section{Derivation of the PLMS algorithm}\label{PLMS_pr}
The objective is to minimize the criterion:
\begin{equation} \min_{\hat{\theta}(t)} J(t) = \textbf{E}\lbrace[e(t)^2]\rbrace
\label{eq:PLMScrit} 
\end{equation}
In order to apply the gradient rule one has to estimate the gradient of the criterion with respect to $\theta (n)$.
\begin{equation}
\bigtriangledown_{\theta} J(t)= 2\textbf{E}\left\{ \frac{\partial e(t)}
{\partial \hat{\theta}(t)}\cdot e(t) \right\} 
\end{equation}
but from \eqref{eq_aposteriori} it results that
\begin{equation}
\frac{\partial e(t)}{\partial \hat{\theta}(t)}=-\phi(t)
\end{equation}
and one obtains the algorithm:
\begin{equation}
 \hat{\theta}(t)=\hat{\theta}(t-1)+\mu \textbf{E}\lbrace\phi(t) e(t)\rbrace
 \label{eq:PLMS1}
\end{equation}
$\textbf{E}\lbrace\phi(t)e(t)\rbrace$ can be approximated by
\begin{equation}
\textbf{E}\lbrace\phi(t) e(t)\rbrace= \frac{1}{N}\sum_{i=0}^{N-1}\phi(t-i)e(t-i)
\end{equation}
and taking $N=1$, the algorithm of \eqref{eq:PLMS1} becomes\footnote{Alternatively one can consider the following approximation:
\begin{equation}
\textbf{E}\lbrace\phi(t) e(t)\rbrace= \frac{1}{N+1}\sum_{i=-N/2}^{N/2}\phi(t+i)e(t+i) \nonumber
\end{equation}
 and one takes $N=0$.}:
\begin{equation}
\hat{\theta}(t)=\hat{\theta}(t-1)+\mu \phi(t) e(t); ~~\mu>0
\label{eq:PLMS2}
\end{equation}
We have to give now an implementable expression of this algorithm. Observe that \eqref{eq_aposteriori} can be re-written as (using \eqref{eq:PLMS2})
\begin{equation}
e(t)= [\theta - \hat{\theta}(t)]^T\phi(t)= e^\circ-\mu\phi(t)^T\phi(t) e(t)
\end{equation}
yielding
\begin{equation}
e(t)=\frac{e^\circ}{1+\mu\phi(t)^T\phi(t)}
\end{equation}
and the adaptation algorithm takes the form given in \eqref{eq_plms} and \eqref{eq_mu}.

\end{appendices}

\bibliographystyle{IEEEtranS}        %
\bibliography{Bibliography}      %

\end{document}